\documentclass[11pt,a4paper]{article}
\usepackage{epsf,epsfig,amsfonts,amsgen,amsmath,amstext,amsbsy,amsopn,amsthm,lineno}
\usepackage{color}
\usepackage{graphicx}
\newtheorem{theorem}{Theorem}

\newtheorem{lemma}{Lemma}

\theoremstyle{definition}

\newtheorem{case}{Case}

\begin{document}
\title
{\bf\Large Degree and neighborhood intersection conditions restricted to induced subgraphs ensuring Hamiltonicity of graphs}

\author{
Bo Ning$^{a}$\thanks{Supported by NSFC (No.~11271300) and the Doctorate Foundation of Northwestern Polytechnical University (No. cx201326). E-mail address: ningbo\_math84@mail.nwpu.edu.cn (B. Ning).}, Shenggui Zhang$^{a}$\thanks{Corresponding author. Supported by NSFC (No.~11271300). E-mail address: sgzhang@nwpu.edu.cn (S. Zhang).} and Bing Chen$^b$\thanks{Supported by NSFC (No.~11271300) and the Scientific Research Program of Shaanxi Provincial Education Department (No. 2013JK0580). E-mail address: cbing\_2004@163.com (B. Chen).}\\
\small $^{a}$Department of Applied Mathematics, School of Science,\\
\small  Northwestern Polytechnical University, Xi'an, Shaanxi 710072, P.R.~China\\[2mm]
\small $^{b}$Department of Applied Mathematics, School of Science, \\
\small Xi'an University of Technology, Xi'an, Shaanxi 710048, P.R.~China}
\date{}
\maketitle

\begin{abstract}
Let claw be the graph $K_{1,3}$. A graph $G$ on $n\geq 3$ vertices is called \emph{o}-heavy if each induced claw of $G$ has a pair of end-vertices with degree sum at least $n$, and called 1-heavy if at least one end-vertex of each induced claw of $G$ has degree at least $n/2$. In this note we show that every 2-connected $o$-heavy or 3-connected 1-heavy graph is Hamiltonian if we restrict Fan-type degree condition or neighborhood intersection condition to certain pairs of vertices in some small induced subgraphs of the graph. Our results improve or extend previous results of Broersma et al., Chen et al., Fan, Goodman \& Hedetniemi, Gould \& Jacobson and Shi on the existence of Hamilton cycles in graphs.
\medskip

\noindent {\bf Keywords:} Hamilton cycle; induced subgraph; claw-free (1-heavy, 2-heavy, \emph{o}-heavy) graph
\smallskip

\noindent {\bf AMS Subject Classification (2000):} 05C38, 05C45
\end{abstract}

\section{Introduction}

We use \cite{Bondy_Murty} for terminology and notation not defined here and consider finite simple graphs only.

Let $G$ be a graph, $H$ be a subgraph of $G$, and $S$ be a subset of $V(G)$. If  $H$ is induced by $S$, then we denote $H$ by $G[S]$. We use $G-H$ to denote the subgraph induced by $V(G)\backslash S$. For a family $\mathcal{H}$ of graphs, $G$ is called $\mathcal{H}$-\emph{free} if $G$ contains no induced subgraph isomorphic to any $H\in\mathcal{H}$. In particular, if $\mathcal{H}=\{H\}$, we use $H$-free instead of $\{H\}$-free.

The graph isomorphic to $K_{1,3}$ is called a \emph{claw}. Its vertex of degree 3 is the \emph{center} and the other vertices are \emph{end-vertices}. In this note, we use claw-free instead of $K_{1,3}$-free. A vertex $v$ of a graph $G$ on $n$ vertices is called \emph{heavy} if $d(v)\geq n/2$. Following \cite{Broersma_Ryjacek_Schiermeyer,Cada}, we say that an induced copy of claw in $G$ is \emph{1-heavy} (\emph{2-heavy}) if at least one (two) of its end-vertices is (are) heavy. The graph $G$ is called \emph{1-heavy} (\emph{2-heavy}) if every induced claw of it is 1-heavy (2-heavy), and called \emph{o-heavy} if every induced claw of it has two end-vertices with degree sum at least $n$. Note that every claw-free graph is 2-heavy, every 2-heavy graph is \emph{o}-heavy, and every \emph{o}-heavy graph is 1-heavy, but the converse does not hold. We use $Z_i$ to denote the graph obtained by attaching a path of length $i$ to one vertex of a triangle.

Generally speaking, there are two types of sufficient conditions with respect to long cycle problems. On the one hand, there are  so-called numerical conditions, among which maybe degree conditions are the most important ones. The following degree condition for Hamilton cycles is due to Fan \cite{Fan}.

\begin{theorem}{\rm (Fan \cite{Fan})}\label{th1}
Let $G$ be a 2-connected graph on $n\geq 3$ vertices. If $\max\{d(u),d(v)\}\geq n/2$ for every pair of vertices $u,v$ in $G$ with $d(u,v)=2$, then $G$ is Hamiltonian.
\end{theorem}

On the other hand, there are so-called structural conditions, of which the forbidden subgraph condition is a good example. The following two results belong to this type of condition.

\begin{theorem}{\rm (Goodman and Hedetniemi \cite{Goodman_Hedetniemi})}\label{th2}
Let $G$ be a 2-connected graph. If $G$ is $\{K_{1,3}, Z_1\}$-free, then $G$ is Hamiltonian.
\end{theorem}

\begin{theorem}{\rm (Gould and Jacobson \cite{Gould_Jacobson})}\label{th3}
Let $G$ be a 2-connected graph. If $G$ is $\{K_{1,3},Z_2\}$-free, then $G$ is Hamiltonian.
\end{theorem}

We also consider the following result due to Shi \cite{Shi} in Hamiltonian graph theory.

\begin{theorem}{\rm (Shi \cite{Shi})}\label{th4}
Let $G$ be a 2-connected graph. If $G$ is claw-free and $|N(u)\cap N(v)|\geq 2$ for every pair of vertices $u$ and $v$ in $G$ with $d(u,v)=2$, then $G$ is Hamiltonian.
\end{theorem}

Let $G$ be a graph on $n\geq 3$ vertices. In the following, we say that a pair of vertices $u,v$ of $G$ with $d(u,v)=2$ satisfies Fan's condition if $\max\{d(u),d(v)\}\geq n/2$, and satisfies Shi's condition if $|N(u)\cap N(v)|\geq 2$.

By restricting Fan's condition and Shi's condition to certain substructures in graphs, Broersma et al. \cite{Broersma_Ryjacek_Schiermeyer} proved the following first result, which generalizes Theorems \ref{th1} and \ref{th4}. Broersma et al. \cite{Broersma_Ryjacek_Schiermeyer} also constructed graphs to show that the condition 2-heavy cannot be replaced by 1-heavy under the same connectedness, and obtained a similar result for 3-connected 1-heavy graphs.

\begin{theorem}{\rm (Broersma, Ryj\'{a}\v{c}ek and Schiermeyer \cite{Broersma_Ryjacek_Schiermeyer})}\label{th5}
Let $G$ be a 2-connected 2-heavy graph. If every pair of vertices $u,v$ with $d(u,v)=2$ satisfies Fan's condition or Shi's condition, then $G$ is Hamiltonian.
\end{theorem}

\begin{theorem}{\rm (Broersma, Ryj\'{a}\v{c}ek and Schiermeyer \cite{Broersma_Ryjacek_Schiermeyer})}\label{th6}
Let $G$ be a 3-connected 1-heavy graph. If every pair of vertices $u,v$ with $d(u,v)=2$ satisfies Fan's condition or Shi's condition, then $G$ is Hamiltonian.
\end{theorem}

Later, Chen et al. \cite{Chen_Zhang_Qiao} extended Theorem \ref{th5} to a larger class of \emph{o}-heavy graphs.

\begin{theorem}{\rm (Chen, Zhang and Qiao \cite{Chen_Zhang_Qiao})}\label{th7}
Let $G$ be a 2-connected $o$-heavy graph. If every pair of vertices $u,v$ in an induced copy $H$ of $Z_1$ with $d_{H}(u,v)=2$ satisfies Fan's condition or Shi's condition, then $G$ is Hamiltonian.
\end{theorem}

Recently, by restricting the degree sum condition to induced claws and the Fan-type degree condition to induced $Z_2$'s, Ning and Zhang \cite{Ning_Zhang} obtained the following theorem.

\begin{theorem}{\rm (Ning and Zhang \cite{Ning_Zhang})}\label{th8}
Let $G$ be a 2-connected $o$-heavy graph. If every pair of vertices $u,v$ in an induced copy $H$ of $Z_2$ with $d_{H}(u,v)=2$ satisfies Fan's condition, then $G$ is Hamiltonian.
\end{theorem}

Motivated by Theorems \ref{th3}, \ref{th7} and \ref{th8}, in this note we first give the following result which is a common generalization of Theorems \ref{th1}, \ref{th3} and \ref{th4}.

\begin{theorem}\label{th9}
Let $G$ be a 2-connected $o$-heavy graph. If every pair of vertices $u,v$ in an induced copy $H$ of $Z_2$ with $d_{H}(u,v)=2$ satisfies Fan's condition or Shi's condition, then $G$ is Hamiltonian.
\end{theorem}

In fact, instead of Theorem \ref{th9}, we prove the following stronger result. Obviously, Theorem \ref{th10} also improves both Theorems \ref{th7} and \ref{th8}.

\begin{theorem}\label{th10}
Let $G$ be a 2-connected $o$-heavy graph. If every pair of vertices $\{a_2,b_1\}$ and $\{a_3,b_1\}$ in an induced $Z_2$ of $G$ satisfies Fan's condition or Shi's condition, then $G$ is Hamiltonian. (see Fig. 1)
\end{theorem}

\begin{center}
\begin{picture}(100,120)
\put(40,30){\circle*{4}} \put(80,30){\circle*{4}} \put(60,60){\circle*{4}} \put(60,90){\circle*{4}} \put(60,120){\circle*{4}}
\put(40,30){\line(1,0){40}} \put(40,30){\line(2,3){20}} \put(80,30){\line(-2,3){20}} \put(60,60){\line(0,1){30}} \put(60,90){\line(0,1){30}} \put(28,30){$a_2$} \put(83,30){$a_3$} \put(65,60){$a_1$} \put(65,90){$b_1$} \put(65,120){$c_1$}
\end{picture}

\small ~~~~~Fig. 1 ~~~The graph $Z_2$.
\end{center}

We also prove the following result which extends Theorem \ref{th6}.

\begin{theorem}\label{th11}
Let $G$ be a 3-connected 1-heavy graph. If every pair of vertices $u,v$ in an induced copy $H$ of $Z_1$ with $d_{H}(u,v)=2$ satisfies Fan's condition or Shi's condition, then $G$ is Hamiltonian.
\end{theorem}

The proofs of Theorems \ref{th10} and \ref{th11} are postponed to Sec. 2.

\section{Proofs of Theorems 10 and 11}
Before the proofs, we give some additional notation and terminology, following \cite{Chen_Zhang_Qiao,Li_Ryjacek_Wang_Zhang}.

Following \cite{Chen_Zhang_Qiao}, let $H$ be a path or a cycle with a given orientation. By $\overleftarrow{H}$ we denote the same graph as $H$ but with the reverse orientation. For a vertex $v$ of $H$, we use $v_H^+$ and $v_H^-$ to denote the immediate successor (if it exists) and immediate predecessor (if it exists) of $v$ on $H$, respectively. Furthermore, we define $v^{+2}_H=(v_H^+)^+$ and $v^{-2}_H=(v_H^-)^{-}$. If $S$ is a set of vertices of $H$, then we define $S^+_H=\{s^+_H|s\in S\}$ and $S^-_H=\{s^-_H|s\in S\}$. When no confusion occurs, we denote $v_H^+$, $v_H^-$, $v_H^{+2}$, $v_H^{-2}$, $S_H^+$ and $S_H^-$ by $v^+$, $v^-$, $v^{+2}$, $v^{-2}$, $S^+$ and $S^-$, respectively. For two vertices $u$ and $v$ of $H$, we use $H[u,v]$ to denote the segment of $H$ from $u$ to $v$ along the orientation.

Let $G$ be a graph on $n\geq 3$ vertices. Recall that a vertex of $G$ is \emph{heavy} if its degree is at least $n/2$; otherwise, it is \emph{light}. For an induced copy of claw in $G$, if all of its three end-vertices are light, it is called a \emph{light claw}. A cycle $C$ in $G$ is called \emph{heavy} if it contains all the heavy vertices of $G$.

Let $G$ be a graph, $C$ a non-Hamilton cycle of $G$ and $v$ a vertex of $V(G)\backslash V(C)$. A subgraph $F$ of $G$ is called a \emph{$(v,C)$-fan} if the following two conditions hold:
(1) $F$ consists of at least two $v$-paths which are mutually vertex-disjoint except for $v$;
(2) the other end-vertices of these $v$-paths are on $C$.
In particular, if $F$ consists of $v$-paths $Q_1,Q_2,\ldots,Q_k$, then we denote $F$ by $(v;Q_1,\ldots,Q_k)$. Thus, if $G$ is $k$-connected, where $k\geq 2$, then by the connectedness condition, there exists a $(v,C)$-fan $F=(v;Q_1,\ldots,Q_k)$.

Next we will introduce a new concept (so called "Ore-cycle" \cite{Li_Ryjacek_Wang_Zhang}) given recently. This concept is motivated by \cite[Lemma 3]{Chen_Zhang_Qiao} and will be used as a main tool in our note. Let $G$ be a graph on $n\geq 3$ vertices and $k\geq 3$ be an integer. Following \cite{Li_Ryjacek_Wang_Zhang}, we use $\widetilde{E}(G)$ to denote the set $\{xy: xy\in E(G)$ or $d(x)+d(y)\geq n, x, y\in V(G)\}$, and call a sequence of vertices $C'=v_1v_2\ldots v_kv_1$ an \emph{Ore-cycle} or briefly, \emph{o-cycle} of $G$, if $v_iv_{i+1}\in \widetilde{E}(G)$ for every $i\in [1,k]$, where the subscripts are taken modulo $k$.

A graph $G$ on $n$ vertices is called \emph{$K_{1,4}$-o-heavy}, if for every induced subgraph $H$ of $G$ isomorphic to $K_{1,4}$, there exist two nonadjacent vertices $u,v\in V(H)$ with degree sum at least $n$.

Three useful lemmas will also be presented. In particular, the last lemma may be a well known fact and appeared in many papers (although not presented as a lemma). For a similar proof, we refer the reader to, for example, \cite{Li_Ryjacek_Wang_Zhang}.
\begin{lemma}[Li \cite{Li}]\label{le1}
Let $G$ be a 2-connected $K_{1,4}$-o-heavy graph. Then every longest cycle is a heavy cycle.
\end{lemma}

\begin{lemma}[Li, Ryj\'{a}\v{c}ek, Wang and Zhang \cite{Li_Ryjacek_Wang_Zhang}]\label{le2}
Let $G$ be a graph on $n$ vertices and $C$ be an o-cycle of $G$. Then there exists a cycle $C'$ of $G$ such that $V(C)\subseteq V(C')$.
\end{lemma}

\begin{lemma}\label{le3}
Let $G$ be a non-Hamiltonian graph, $C$ be a longest cycle (longest heavy cycle) of $G$, $R$ a component of $G-V(C)$, and $A=\{v_1,v_2,\ldots,v_k\}$ the set of neighbors
of $R$ on $C$. For $u\in R$, $v_i,v_j\in A$, there hold\\
$(a)$ $uv_i^-\notin \widetilde{E}(G)$ and $uv_i^+\notin \widetilde{E}(G)$;\\
$(b)$ $v_i^-v_j^-\notin \widetilde{E}(G)$ and $v_i^+v_j^+\notin \widetilde{E}(G)$; and\\
$(c)$ $v_iv_j^-\notin \widetilde{E}(G)$ and $v_iv_j^+\notin \widetilde{E}(G)$ if $v^-_iv_i^+\in \widetilde{E}(G)$.\\
Furthermore, if $G$ is a 2-connected $o$-heavy graph, then there hold\\
$(d)$ $v_i^{-}v_i^{+}\in \widetilde{E}(G)$ and $v_j^{-}v_j^{+}\in \widetilde{E}(G)$;\\
$(e)$ $v^-_iv_i^+\in E(G)$ or $v^-_jv_j^+\in E(G)$; and\\
$(f)$ Let $y_1$ be the first vertex on $\overrightarrow{C}[v_1,v_2]$ such that $v_1y_1\notin E(G)$ and $y_2$ be the first vertex on $\overrightarrow{C}[v_2,v_1]$ such that $v_2y_2\notin E(G)$. For $u\in V(P)\backslash \{v_1,v_2\}$, $w_1\in \overrightarrow{C}[v_1^+,y_1]$ and $w_2\in \overrightarrow{C}[v_2^+,y_2]$, we have
(1) $uw_1,uw_2\notin \widetilde{E}(G)$; (2) $v_1w_2,v_2w_1\notin \widetilde{E}(G)$; and (3) $w_1w_2\notin \widetilde{E}(G)$.
\end{lemma}

\noindent{}
{\bf Proof of Theorem \ref{th10}.}
We prove the theorem by contradiction. Suppose that $G$ satisfies the condition of Theorem \ref{th10} and is not Hamiltonian. Let $C$ be a longest cycle of $G$. Since every $o$-heavy graph is also $K_{1,4}$-$o$-heavy, by Lemma \ref{le1}, $C$ is a heavy cycle. Since $G$ is not Hamiltonian, $G-C\neq \emptyset$ and every vertex in $G-C$ is light. Let $H$ be a component of $G-C$. By the connectedness condition, there is a $(v_1,v_2)$-path of length at least 2 with all internal vertices in $H$, where $v_1,v_2\in V(C)$. We choose such a path $P$ as short as possible. Let $P=u_0u_1u_2\ldots u_ru_{r+1}$, where $u_0=v_1,u_{r+1}=v_2$.

We divide the proof into two cases.

\begin{case}
$r\geq 2$, or $r=1$ and $v_1v_2\notin E(G)$.
\end{case}
By Lemma \ref{le3} ($b$), we obtain $d(v_1^-)+d(v_2^-)<n$ and $d(v_1^+)+d(v_2^+)<n$. This implies that $d(v_1^-)+d(v_1^+)<n$ or $d(v_2^-)+d(v_2^+)<n$. W.l.o.g., assume that $d(v_1^-)+d(v_1^+)<n$ and $d(v_1^-)<n/2$. By Lemma \ref{le3} ($a$), $u_1v_1^-\notin E(G)$ and $u_1v_1^+\notin E(G)$. Since $G$ is $o$-heavy, we obtain $v_1^-v_1^+\in E(G)$. If $u_2\neq v_2$, then by Lemma \ref{le3} ($a$), there are no edges between $\{u_1,u_2\}$ and $\{v_1^-,v_1^+\}$, and by the choice of $P$, $u_2v_1\notin E(G)$. If $u_2=v_2$, then by the condition, $v_1u_2\notin E(G)$, and it follows from Lemma \ref{le3} ($c$) that $u_2v_1^-,u_2v_1^+\notin E(G)$. Thus the subgraph induced by $\{v_1,v_1^-,v_1^+,u_1,u_2\}$ is a $Z_2$. Noting that both $v_1^-$ and $u_1$ are light, $v_1^-$ and $u_1$ have a common neighbor $z$ other than $v_1$. By Lemma \ref{le3} ($a$), $v_1^-\notin N_C(H)$, and this implies that $z\in V(C)$. Furthermore, $u_1z^-\notin \widetilde{E}(G)$ and $u_1z^+\notin \widetilde{E}(G)$ by Lemma \ref{le3} ($a$). Since $G$ is $o$-heavy, $z^-z^+\in \widetilde{E}(G)$, and in this case, the fact $v_1^-z\in E(G)$ contradicts to Lemma \ref{le3} ($c$).

\begin{case}
$r=1$ and $v_1v_2\in E(G)$.
\end{case}
Let $y_1$ be the first vertex in $\overrightarrow{C}[v_1^+,v_2^-]$ that is nonadjacent to $v_1$, and $y_2$ be the first vertex in $\overrightarrow{C}[v_2^+,v_1^-]$ that is nonadjacent to $v_2$. Note that $v_1v_2^-\notin E(G)$ and $v_2v_1^-\notin E(G)$. Thus, $y_1,y_2$ are well-defined.

By Lemma \ref{le3} ($f$), $d(y_1^-)+d(y_2^-)<n$, and this implies either $y_1^-$ or $y_2^-$ is light. W.l.o.g, assume that $y_1^-$ is light. By Lemma \ref{le3} ($f$), there are no edges between $\{v_2,u_1\}$ and $\{y_1^-,y_1\}$. Thus the subgraph induced by $\{v_1,u_1,v_2,y_1^-,y_1\}$ is a $Z_2$. Since both $u_1$ and $y_1^-$ are light, $u_1$ and $y_1^-$ have a common neighbor $z$ other than $v_1$.

By Lemma \ref{le3} ($f$), we have $z\notin V(H)$ and $z\notin \overrightarrow{C}[v_1^+,y_1]$. This implies that $z\in \overrightarrow{C}[y_1^+,v_1^{-2}]$. By Lemma \ref{le3} ($a$), $u_1z^-\notin E(G)$ and $u_1z^+\notin E(G)$. Since $G$ is $o$-heavy, $z^-z^+\in \widetilde{E}(G)$. If $y_1^-=v_1^+$, then
$C'=\overrightarrow{C}[v_1^+,z^-]z^-z^+\overrightarrow{C}[z^+,v_1]v_1u_1zv_1^+$ is an $o$-cycle longer than $C$. If $y_1^-\neq v_1^+$, then $v_1y_1^{-2}\in E(G)$, and $C'=\overrightarrow{C}[y_1^-,z^-]z^-z^+\overrightarrow{C}[z^+,v_1^-]v_1^-v_1^+\overrightarrow{C}[v_1^+,y_1^{-2}]$
$y_1^{-2}v_1u_1zy_1^-$ is an $o$-cycle longer than $C$. In each case, we can get an $o$-cycle longer than $C$ and contains all vertices in $C$. By Lemma \ref{le2}, there is a cycle longer than $C$, contradicting the choice of $C$.

The proof is complete.
{\hfill$\Box$}

\noindent{}
{\bf Proof of Theorem \ref{th11}.} We prove the theorem by contradiction. Suppose that $G$ satisfies the condition of Theorem \ref{th11} and is not Hamiltonian. We choose $C$ as a longest cycle in $G$ with a given orientation. Since $G$ is not Hamiltonian, $G-C\neq \emptyset$. Obviously, every 1-heavy graph is also $K_{1,4}$-$o$-heavy, then by Lemma \ref{le1}, $C$ is a heavy cycle and this implies that every vertex in $G-C$ is light. Let $H$ be a component of $G-C$. Now we choose $F=(u;Q_1,Q_2,Q_3)$ as a $(u,C)$-fan satisfying $|V(F)|$ is as small as possible, where $u$ is chosen over all the vertices in $H$. Let $v_i$ be the end-vertex of $Q_i$ on $C$ along the orientation of $C$, respectively. Set $u_i={v_i}_{Q_i[u,v_i]}^-$, where $i=1,2,3$.

From Lemma \ref{le3} ($b$), there is at most one heavy vertex in $N_{C}^+(H)$ and at most one heavy vertex in $N_{C}^-(H)$. Thus, we can assume $v_1^{-},v_1^{+}$ are light. By Lemma \ref{le3} ($a$), $u_1v_1^-\notin E(G)$ and $u_1v_1^+\notin E(G)$. If $v_1^-v_1^+\notin E(G)$, then $\{v_1,v_1^-,v_1^+,u_1\}$ induces a claw. Note that $v_1^{-},v_1^{+}$ are light and $u_1$ is light by the choice of $C$ and Lemma \ref{le1}. Thus $\{v_1,v_1^-,v_1^+,u_1\}$ induces a light claw, contradicting the fact that $G$ is 1-heavy. So, $v_1^-v_1^+\in E(G)$.

Note that $\{v_1,v_1^-,v_1^+,u_1\}$ induces a $Z_1$.  Since $v_1^-,v_1^+,u_1$ are light, there exist vertices $z_1,z_2$ such that $z_1\in N(u_1)\cap N(v_1^-)\backslash \{v_1\}$ and $z_2\in N(u_1)\cap N(v_1^+)\backslash \{\textcolor{red}{v_1}\}$. By Lemma \ref{le3} ($a$), $z_i\in V(C)$, where $i=1,2$.

Furthermore, we claim that $z_1\neq z_2$. If not, then $\{z_1,u_1,v_1^-,z_1^-\}$ and $\{z_1,u_1,v_1^+,z_1^+\}$ induce two claws. Since $G$ is 1-heavy, $z_1^-,z_1^+$ are heavy, and this implies  $z_1^-z_1^+\in \widetilde{E}(G)$. However, $v_1^-z_1\in E(G)$, contradicting Lemma \ref{le3} ($c$). Thus, $z_1\neq z_2$ and we obtain $z_1^-,z_2^+$ are heavy.

Now we can choose $F=(u_1;u_1v_1,u_1z_1,u_1z_2)$ (where $v_1,z_1,z_2$ appear on $C$ along its orientation) or $F=(u_1;u_1v_1,u_1z_2,u_1z_1)$ (where $v_1,z_2,z_1$ appear on $C$ along its orientation) as the required fan.

If $F=(u_1;u_1v_1,u_1z_1,u_1z_2)$ is the required fan, then noting that $z_1^-z_2^+\in \widetilde{E}(G)$ and $v_1^-z_1,v_1^+z_2\in E(G)$, $C'=u_1z_1\overrightarrow{C}[z_1,z_2]z_2v_1^+\overrightarrow{C}[v_1^+,z_1^-]z_1^-z_2^+\overrightarrow{C}[z_2^+,v_1]v_1u_1$ is an $o$-cycle including all vertices in $C$ and longer than $C$. By Lemma \ref{le2}, there is a cycle longer than $C$, a contradiction.

If $F=(u_1;u_1v_1,u_1z_2,u_1z_1)$ is the required fan, then $z_2^+,z_1^-$ are heavy and this implies $z_2^-,z_1^+$ are light by Lemma \ref{le3} ($b$). If $v_1^-z_1^+\notin E(G)$, then $\{z_1,u_1,z_1^+,v_1^-\}$ induces a light claw, a contradiction. Thus, $v_1^-z_1^+\in E(G)$. By Lemma \ref{le3} ($b$), ($c$), $\{v_1^-,v_1,v_1^+,z_1^+\}$ induces a $Z_1$. Since $v_1^+,z_1^+$ are light, there is a vertex other than $v_1^-$, say $w$, is a common neighbor of $z_1^+$ and $v_1^+$. By Lemma 3 ($a$), $w\in V(C)$. If $w=z_1$ or $w=z_2$, then $\{w,v_1^+,z_1^+,u_1\}$ induces a light claw, a contradiction. Thus, $w\neq z_1,z_2$. Consider the subgraph induced by $\{w,w^+,z_1^+,v_1^+\}$. Since $G$ is 1-heavy and $z_1^+,v_1^+$ are light, $w^+$ is heavy. Thus, $w^+z_1^+\in \widetilde{E}(G)$, a contradiction to Lemma \ref{le3} ($b$).

The proof is complete. {\hfill$\Box$}


\end{document}